\documentclass[12pt]{article}
\usepackage{amsmath}
\usepackage{amssymb}
\usepackage{amsthm}

\usepackage[margin=2.0cm]{geometry}
\usepackage{hyperref}
\usepackage{multirow}

\usepackage{color}
\usepackage{enumerate}
\usepackage{tikz,lmodern}
\usepackage{amsfonts}
\usetikzlibrary{arrows, decorations.markings}

  \def\ort{\overrightarrow}  
\def\Rm#1{\lowercase\expandafter{\romannumeral#1}}

  \def\lg{\langle} \def\rg{\rangle}\def\Ga{\Gamma}
\def\Hol{\hbox{\rm Hol}}
  \def\D{\Delta} 
  
\def\a{\alpha} \def\b{\beta} \def\g{\gamma} \def\d{\delta} 
 \def\s{\sigma} \def\t{\tau}  
     
 \def\Aut{\hbox{\rm Aut}}   
    \def\D{\hbox{\rm D}} 
  \def\mod{\hbox{\rm mod}\,}

\def\Cay{\hbox{\rm Cay}}  
  
\def\qed{\hfill $\Box$}

\def\Hol{\hbox{\rm Hol}} 
\newtheorem{theorem}{Theorem}[section]%
\newtheorem{lemma}[theorem]{Lemma}%
\newtheorem{proposition}[theorem]{Proposition}%
\def\mz{{\mathbb Z}}

\begin{document}
\title{Normal and non-normal Cayley digraphs on cyclic and dihedral groups}
\author{\\ Jun-Feng Yang, Yan-Quan Feng\footnotemark, Fu-Gang Yin, Jin-Xin Zhou\\
{\small\em School of mathematics and statistics, Beijing Jiaotong University, Beijing, 100044, P.R. China}\\
}

\renewcommand{\thefootnote}{\fnsymbol{footnote}}
\footnotetext[1]{Corresponding author.
E-mails:
jf.yang@bjtu.edu.cn (J.-F. Yang),
yqfeng@bjtu.edu.cn (Y.-Q. Feng),
fgyin@bjtu.edu.cn (F.-G. Yin),
jxzhou@bjtu.edu.cn (J.-X. Zhou)}

\date{}
\maketitle

\begin{abstract}
  A Cayley digraph on a group $G$ is called NNN if the Cayley digraph is normal and its automorphism group contains a non-normal regular subgroup isomorphic to $G$. A group is called NNND-group or NNN-group if there is an NNN Cayley digraph or graph on the group, respectively. In this paper, it is shown that there is no cyclic NNND-group, and hence no cyclic NNN-group. Furthermore, a dihedral group of order $2n$ is an NNND-group or an NNN-group if and only if $n\ge 6$ is even and $n\not=8$.
  
\end{abstract}

\section{Introduction}
Throughout this paper, digraphs are finite and simple, and groups are finite.
For a digraph $\Gamma$, $\Aut(\Gamma)$ denotes the automorphism group of $\Gamma$, and $V(\Gamma)$ and $A(\Gamma)$ denote the vertex set and arc set of $\Gamma$ respectively, while $V(\Gamma)$ and $E(\Gamma)$ denote the vertex set and edge set of $\Gamma$ respectively when $\Gamma$ is a graph. For convenience, we view a graph $\Gamma$ with vertex set $V(\Gamma)$ and edge set $E(\Gamma)$ as a digraph with vertex set $V(\Gamma)$ and arc set $\{(u,v)\ |\ \{u,v\}\in E(\Gamma)\}$, and this makes no confusion in this paper. Thus, a graph is a special case of a digraph.

Let $G$ be a finite group and let $S$ be a subset of $G$ with $1\notin S$. The {\em Cayley digraph $\Cay(G,S)$} on $G$ with respect to $S$ is defined to be the digraph with vertex set $G$ and arc set $\{(g,sg)\ |\ g\in G, s\in S\}$. Let $\Gamma=\Cay(G, S)$. It is easy to see that $\Gamma$ is connected if and only if $S$ generates $G$. When $S$ is inverse closed, that is $S=S^{-1}$, $\Gamma$ is called a {\em Cayley graph}, and in this case, $\Gamma$ has edge set $\{\{g,sg\}\ |\ g\in G, s\in S\}$. % and denote by $\Aut(\Cay(G,S))$ the automorphism group of $\Cay(G,S)$.
For a given $g\in G$, define the map $R(g):x\mapsto xg$ for each $x\in G$, and write $R(G)=\{R(g)\ |\ g\in G\}$, called {\em the right regular representation of $G$}. Then $R(G)$ is a regular subgroup of $\Aut(\Gamma)$, that is, $R(G)$ is vertex transitive on $G$ and only identity can fix an element in $G$. Moreover, a digraph is isomorphic to a Cayley digraph $\Sigma$ on a group $G$ if and only if $\Aut(\Sigma)$ has a regular subgroup that is isomorphic to $G$ (see \cite{DMM}).

Determining the automorphism groups of Cayley graphs is one of the central issues in the graph theory. Godsil~\cite{Godsil} proved that $R(G)\rtimes\Aut(G,S)$ is just the normalizer of $R(G)$ in $\Aut(\Cay(G,S))$, where $\Aut(G, S)=\{\alpha\in\Aut(G)\ |\ S^\alpha=S\}$. This plays an important role in the study of the automorphisms of Cayley graphs. Based on this, Xu~\cite{XuM} introduced the so-called normal Cayley digraph. A Cayley digraph $\Cay(G,S)$ is called {\em normal} if $R(G)$ is a normal subgroup of $\Aut(\Cay(G, S))$. %By Godsil~\cite{Godsil}, if $\Cay(G,S)$ is normal, then $\Aut(\Cay(G,S))=R(G)\rtimes \Aut(G,S)$.

The isomorphism problem for Cayley graphs is also a fundamental problem in graph theory. A Cayley digraph $\Cay(G,S)$ is called {\em CI} (CI stands for Cayley isomorphism) if, for any Cayley digraph $\Cay(G,T)$, whenever $\Cay(G,S)\cong \Cay(G,T)$, we have $S^{\s}=T$ for some $\s\in\Aut(G)$. There has been quite a lot of research on CI-graphs, for which we refer the reader to the nice survey paper~\cite{Li}. By Babai~\cite{Ba77}, a Cayley digraph $\Gamma=\Cay(G,S)$ on a group $G$ is CI if and only if all regular subgroups of $\Aut(\Gamma)$ isomorphic to $G$ are conjugate. From this we know that if $\Gamma$ is a normal Cayley digraph on $G$, then $\Gamma$ is CI if and only if $R(G)$ is the unique regular subgroup of $\Aut(\Gamma)$ isomorphic to $G$ (see \cite[Corollary~6.9]{Li}). Li~\cite{Li} gave some examples of normal Cayley digraphs which are non-CI. However, he commented that ``normal Cayley graphs which are not CI-graphs seem to be very rare" and he proposed the problem: {\em Characterize normal Cayley digraphs which are not CI-graphs.}

The first example of normal Cayley graph on elementary abelian group which is non-CI was found by Royle~\cite{Royle} in 2008. Royle's graph is a normal Cayley graph on $\mz_2^6$ whose automorphism group also contains a non-normal regular subgroup isomorphic to $\mz_2^6$. Such a Cayley graph is called an NNN graph. More precisely, a Cayley digraph $\Gamma=\Cay(G,S)$ on a group $G$ is said to be {\em NNN} if $\Gamma$ is normal and $\Aut(\Gamma)$ has a non-normal regular subgroup isomorphic to $G$. In 2010, Giudici and Smith \cite{GS} constructed a new NNN graph which is a strongly regular Cayley graph on $\mz_6^2$. The first infinite family of NNN graphs were constructed by Bamberg and Giudici \cite{BM} in 2011 as the point graphs of a particular family of generalized quadrangles. In 2017, Xu~\cite{XuY17} constructed some more infinite families of NNN graphs by using Cartesian product, direct product and strong product.

Another interesting problem in the study of NNN graphs is to characterize the groups which admit NNN digraphs. A group $G$ is called {\em NNND-group} or {\em NNN-group} if there exists an NNN Cayley digraph or graph on $G$, respectively. As we view a graph as a digraph, an NNN-group is an NNND-group. In 2018, Xu proved in \cite{XuY18} that every elementary abelian $2$-group $\mz_2^d$ with $d\geq 6$ is an NNN-group, and that there is no non-abelian simple NNND-group. In 2022,  Xu~\cite{XuY22} showed that the symmetric group $S_n$ of degree $n$ is an NNN-group if and only if $n\ge 5$. In 2019, Giudici, Morgan and Xu~\cite{GLX} proved that every cyclic group is not an NNN-group.

In this paper, we shall give a classification of cyclic and dihedral NNND-groups. In the literature, the isomorphism problem for the Cayley digraphs on cyclic or dihedral groups have been extensively studied (see, for example, \cite{Ba77,conder,DT,Kov,Kovacs,C.H.Li9,Muzychuk, M.Muzychuk,XFRL,XFK,XFZ}). This provides some motivation for us to investigate the cyclic and dihedral NNND-groups.

Our first result shows that there is no cyclic NNND-group.

\begin{theorem}\rm\label{mainth1}
Every cyclic group is not an NNND-group.
\end{theorem}

\noindent{\bf Remark.}\ (1)\ Note that every NNN-group is also an NNND-group. From Theorem~\ref{mainth1} we immediately obtain Giudici-Morgan-Xu's result, namely, every cyclic group is not an NNN-group.

(2)\ The main approach of Giudici, Morgan and Xu~\cite{GLX} is to analyze the regular subgroups in the holomorph of cyclic groups. In this paper, we shall present a totally different and much shorter proof for this result.\medskip

Our second result is the classification of dihedral NNND-groups and NNN-groups.

\begin{theorem}\rm\label{mainth2}
Let $\D_{2n}$ be a dihedral group of order $2n$. Then the following statements are equivalent:
  \begin{itemize}
    \item[(1)] $\D_{2n}$ is an NNN-group;
    \item[(2)] $\D_{2n}$ is an NNND-group;
    \item[(3)] $n\geq 6$ is even, and $n\not=8$.
  \end{itemize}
\end{theorem}

The layout of this paper is as follows. In section 2, we give some definitions and preliminary results for later use. Theorem~\ref{mainth1} is proved in Section~\ref{theorem1} and  Theorem~\ref{mainth2} is proved in Section~\ref{theorem2}.

\section{Preliminaries}

In this section, we present some notation and results that will be used later.
Let $\Ga=\Cay(G,S)$ be a Cayley digraph on a group $G$. Then $\Aut(G,S)$ is a subgroup of $\Aut(\Ga)_1$, the stabilizer of $1$ in $\Aut(\Ga)$. The following proposition is a criterion for normal Cayley digraph.

\begin{proposition}\rm\label{normal}(\cite[Proposition 1.5]{XuM}) Let $\Ga=\Cay(G,S)$ be a Cayley digraph on a group $G$ with respect to $S$ and let $A=\Aut(\Ga)$. Then $\Cay(G,S)$ is normal if and only if $A=R(G)\rtimes \Aut(G,S)$, and if and only if $A_1=\Aut(G,S)$.
\end{proposition}

%By Babai's criterion (see \cite{Ba77}), if a Cayley digraph $\Ga=\Cay(G,S)$ is CI, then all regular subgroups of $\Aut(\Ga)$ isomorphic to $G$ are conjugate in $\Aut(\Ga)$, which implies that a CI Cayley digraph cannot be NNN. %A group $G$ is called {\em DCI-group} or {\em CI-group} if all Cayley digraphs or graphs on $G$ are CI, respectively. It follows that DCI-group and CI-group are not NNND-group and NNN-group, respectively. Furthermore,
A group is called an {\em NDCI-group} or {\em NCI-group}, if all normal Cayley digraphs or graphs on the group are CI, respectively. By definition,  every NDCI-group is not an NNND-group while every NCI-group is not an NNN-group.

The following proposition is a classification of cyclic NDCI-groups and NCI-groups.

\begin{proposition}\rm\label{NDCI}  (\cite[Theorem 1.1]{XFRL})
A cyclic group of order $n$ is an NDCI-group if and only if $8\nmid n$, and is an NCI-group if and only if either $n=8$ or $8\nmid n$.
\end{proposition}

The following proposition is a classification of dihedral NDCI-groups and NCI-groups.

\begin{proposition}\rm\label{NDCID}(\cite[Theorem 1.1]{XFZ})
  Let $n\ge 2$ be an integer and let $\D_{2n}$ be the dihedral group of order $2n$.
  Then the following statements are equivalent:
  \begin{itemize}
    \item[(1)] $\D_{2n}$ is an NDCI-group;
    \item[(2)] $\D_{2n}$ is an NCI-group;
    \item[(3)] Either $n=2,4$ or $n$ is odd.
  \end{itemize}
\end{proposition}

Let $G$ be a group and let $L\subseteq \Aut(G)$.
Write $F_G(L)=\{g\  |\ g^l=g \  {\rm for \ all}\  l\in L\}$, the set of fixed-points of $L$ in $G$. Then $F_G(L)\leq G$.
The following proposition offers some sufficient conditions for a Cayley digraph being non-normal.

\begin{proposition}\rm\label{MainPro} (\cite[Theorem 2.2]{XFK})
  Let $\Cay(G,S)$ be a Cayley digraph. Let $1\neq L\le \Aut(G,S)$ and let $K\unlhd G$ such that for every right coset $Kg$ in $G$, either $L$ fixes $Kg$ pointwise, or $Kg$ is an orbit of $L$. Assume that one of the following holds:
  \begin{itemize}
    \item[(1)] $|G:F_G(L)|>2$;
    \item[(2)] $|G:F_G(L)|=2$, and there are $g\in G\setminus F_G(L)$ and $k\in K$ such that $k^g\neq k^{-1}$;
    \item[(3)] $|G:F_G(L)|=2$, and there is $1\neq \gamma\in \Aut(G,S)$ such that $F_G(\lg \gamma\rg)\neq F_G(L)$ and $\gamma$ fixes every coset of $K$ in $G$ setwise.
  \end{itemize}
  Then $\Cay(G,S)$ is non-normal.
\end{proposition}

%Let $G$ be a group. For a given $g\in G$, define $L(g):x\mapsto g^{-1}x$ for all $x\in G$, and let $L(G)=\{L(g)\ |\ g\in G\}$, called the {\em left regular representation} of $G$. Let $\Ga=\Cay(G,S)$ be a Cayley graph on $G$. Then   $S=S^{-1}$. We know that $R(G)\leq \Aut(\Ga)$, but $L(G)$ may not be a subgroup of $\Aut(\Ga)$. If $L(G)\leq \Aut(\Ga)$, $\Ga$ is called a {\em dual Cayley graph}. The following proposition outlines some essential properties of dual Cayley graph.

%\begin{proposition}\rm\label{DCG}(\cite[Theorem 2.2]{Pan20})
% Let $\Ga=\Cay(G,S)$ be a dual Cayley graph. The following statements are equivalent:
% \begin{itemize}
%   \item[(1)] $L(G)\le \Aut(\Ga)$;
 %  \item[(2)] $\t\in \Aut(\Ga)$ where $\t: g\mapsto g^{-1}$ for all $g\in G$;
  % \item[(3)] $S=S^{-1}$ is a union of some $G$-conjugate classes.
% \end{itemize}
%\end{proposition}

%\noindent {\bf Remark}: Note that a group $G$ is abelian if and only if $\tau$, as given in Proposition~\ref{DCG}~(2), is an automorphism of $G$. By Propositions~\ref{DCG}~(2) and \ref{normal}, a dual Cayley graph on a non-abelian group is non-normal.

\section{Proof of Theorem~\ref{mainth1}\label{theorem1}}

For a positive integer $n$, let
\begin{equation}\label{factorizationofn}
n=\prod_{i=1}^{s+1}{p_i^{k_i}}, \mbox{ with } p_1>p_2>\cdots>p_s>p_{s+1}=2,
\end{equation}
represent the distinct prime factorization of $n$, and let $C_n$ be the cyclic group of order $n$. Then $p_i$ is an odd prime for every $1\leq i\leq s$. For $x\in C_n$, denote by $o(x)$ the order of $x$ in $C_n$. Then we may write
\begin{equation}\label{c_nfactorization}
C_n=C_{p_1^{k_1}}\times C_{p_2^{k_2}}\times\cdots\times C_{p_s^{k_s}}\times C_{p_{s+1}}^{k_{s+1}}=\lg a\rg=\lg a_1\rg \times\cdots \times \lg a_s\rg\times \lg a_{s+1}\rg,
\end{equation}
where $o(a)=n$ and $o(a_i)=p_i^{k_i}$ for every $1\le i\le s+1$.
By \cite[Theorem 4.7]{Rot},
\begin{equation}\label{autoc_nfactorization}
\Aut(C_n)=\Aut(C_{p_1^{k_1}}) \times \Aut(C_{p_2^{k_2}})\times \cdots \times \Aut(C_{p_s^{k_s}})\times \Aut(C_{p_{s+1}^{k_{s+1}}}),
\end{equation}
where $\Aut(C_{p_i^{k_i}})$ is identified as the subgroup of $\Aut(C_n)$ by defining for every $\d_i\in \Aut(C_{p_i^{k_i}})$ that $a_j^{\d_i}=a_j$ for all $j\neq i$.
Furthermore, $\Aut(C_{p_i^{k_i}})\cong C_{p_i^{k_i-1}(p_i-1)}$ for $1\leq i\leq s$ as $p_i$ is an odd prime, and $\Aut(C_{p_i^{k_i}})$ has an element of order $p_i$ if and only if $k_i\geq 2$. Since $p_{s+1}=2$, we have $\Aut(C_{p_{s+1}^{k_{s+1}}})\cong C_2\times C_{2^{k_{s+1}-2}}$ for $k_{s+1}\geq 2$.

Let $k_i\geq 2$ and let $\a_i$ be the automorphism of $C_n$ of order $p_i$  induced by
\begin{equation}\label{auorderp_i}
\a_i: a_i\mapsto a_i^{p_i^{k_i-1}+1}, \mbox{ and } a_j\mapsto a_j \mbox{ for every } j\neq i.
\end{equation}
If $1\leq i\leq s$, then $\Aut(C_{p_i^{k_i}})$ is cyclic as $p_i$ is odd, which implies that $\lg \a_i\rg$ is the unique subgroup of order $p_i$ in $\Aut(C_{p_i^{k_i}})$.

Let $k_{s+1}\geq 4$. Then $\Aut(C_{p_{s+1}^{k_{s+1}}})$ has an element of order $4$. Let $\b$ be the automorphism of $C_n$ of order $4$ induced by
\begin{equation}\label{auorder4}
\b: a_{s+1}\mapsto a_{s+1}^{2^{k_{s+1}-2}+1}, \mbox{ and } a_j\mapsto a_j \mbox{ for every } j\neq s+1.
\end{equation}

It is easy to see that $\b^2=\a_{s+1}$. Throughout this section, we use the notations and formulae in Eqs~(\ref{factorizationofn})-(\ref{auorder4}). First we have the following result.

\begin{lemma}\rm\label{pi}
Let $\Cay(C_n,S)$ be a Cayley digraph on the cyclic group $C_n$. Assume $\a_i\in \Aut(C_n,S)$ for some $1\leq i\leq s$ or $\b\in \Aut(C_n,S)$. Then $\Cay(C_n,S)$ is non-normal.
\end{lemma}

\proof First assume that $\a_i\in \Aut(C_n,S)$ for some $1\leq i\leq s$. Then $k_i\geq 2$. Set $K=\lg a^{n/p_i}\rangle=\lg  a_i^{p_i^{k_i-1}}\rg$  and $H=\lg  a_1,\cdots, a_{i-1},a_{i+1},\cdots, a_s, a_{s+1}\rg$. Then $K$ has order $p_i$ and
$H$ is the Hall $p_i'$-subgroup of $C_n$. Furthermore, $C_n=\cup_{x\in H}\cup_{k=0}^{p_i^{k_i-1}-1}xa_i^kK$, and $K\unlhd C_n$ as $C_n$ is abelian.

Let $\Ga=\Cay(C_n,S)$ and write $L=\lg \a_i \rg$. By assumption,  $L\leq \Aut(C_n,S)\leq \Aut(\Ga)$. Denote by $F_{C_n}(L)$ the fixed point set of $L$ in $C_n$. By Eq~(\ref{auorderp_i}), it is easy to see that $F_{C_n}(L)=\lg a_1,\cdots, a_{i-1}, a_i^{p_i}, a_{i+1},\cdots, a_s\rg$. It follows that $|C_n:F_{C_n}(\lg L \rg)|=p_i>2$ as $p_i$ is odd ($1\leq i\leq s$), and $K\leq F_{C_n}(L)$ as $k_i\geq 2$.  It is easy to see that for every $x\in H$, $L$ fixes $xa_i^kK$ pointwise when $k \equiv 0\ \mod p_i$, and  $xa_i^kK$ is an orbit of $L$ when $k \not\equiv 0\ \mod p_i$. By Proposition~\ref{MainPro}, $\Ga$ is non-normal.

Now assume  $\b\in \Aut(C_n,S)$. Then $k_{s+1}\geq 4$ and $\a_{s+1}=\b^2\in  \Aut(C_n,S)$. Set $E=\langle a_1a_2\cdots a_sa_{s+1}^4\rangle$ and $I=\langle a_{s+1}^{2^{k_{s+1}-1}}\rangle$. Then $E$ and $I$ are the unique subgroups of order $n/4$ and $2$ in $C_n$, respectively. Furthermore, $C_n=E_0\cup E_1\cup E_2\cup E_3$, where $E_i=a_{s+1}^iE$ for $i=0,1,2$ or $3$. Thus, $E_0=E$, and $I\leq E$ as $k_{s+1}\geq 4$. Taking $i=s+1$ in Eq~(\ref{auorderp_i}), it is easy to see that $\a_{s+1}$ fixes $E_0\cup E_2$ pointwise, and interchanges the two elements in every coset of $I$ in $E_1\cup E_3$. By Eq~(\ref{auorder4}), $\b$ fixes every vertex in $E_0$, and interchanges the two elements in every coset of $I$ in $E_2$. For convenience, denote by $4K_1$ the graph with $4$ isolated vertices, and by $\ort{K}_{2,2}$ the complete bipartite digraph of order $4$, which has vertex set $\{1,2,3,4\}$ and arc set $\{(1\ 3), (1\ 4),(2\ 3),(2\ 4)\}$.
For any cosets $xI$ and $yI$ such that $x\in E_i$ and $y\in E_j$ with $i\not=j$, we have $xI\subseteq E_i$ and $yI\subseteq E_j$. Let $(xI,yI)$ be the induced subdigraph of $\Gamma$ from $xI$ to $yI$, that is, $V((xI,yI))=xI\cup yI$ and $A((xI,yI))=\{(u,v)\ |\ u\in xI, v\in yI, (u,v)\in A(\Gamma)\}$.

\medskip
\noindent {\bf Claim:} $(xI,yI)\cong 4K_1$ or $\ort{K}_{2,2}$ for any $x\in E_i$ and $y\in E_j$ with $i\not=j$.

Let $(1,z)\in A(\Gamma)$ for some $z\in C_n\setminus E$. Note that $1\in E_0=E\leq C_n$. Then $z\in E_1\cup E_2\cup E_3$. Assume that $z\in E_1\cup E_3$. Recall that $\a_{s+1}$ interchanges the two elements in $zI$ as $zI\subseteq E_1\cup E_3$. Since $\a_{s+1}\in \Aut(C_n,S)$, there is an arc from $1$ to every vertex in  $zI$, and hence $[I,zI]\cong \ort{K}_{2,2}$ as $R(I)\leq R(C_n)\leq \Aut(\Gamma)$. Since $R(E)\leq \Aut(\Gamma)$ and $R(E)$ fixes $E_i$ setwise for every $0\leq i\leq 3$, we have $(xI,yI)\cong 4K_1$ or $\ort{K}_{2,2}$ for any $x\in E_0$ and $y\in E_j$ with $j=1$ or $3$. This is also true for $j=2$ because $\b$ fixes every vertex in $E_0$ and interchanges the two elements in every coset of $I$ in $E_2$. The claim follows from the fact that $R(a_{s+1})\in \Aut(\Gamma)$ maps $E_i$ to $E_{i+1}$ for every $0\leq i\leq 2$ and maps $E_3$ to $E_0$.    \qed

\medskip

By Claim, if $\gamma\in \Aut(\Gamma)$ fixes every coset of $I$ in $C_n$, then the restriction $\gamma_{E_i}$ of $\gamma$ on $E_i$ can be extended to an automorphism, denoted by $\gamma_i$, of $\Gamma$ by fixing every vertex in $C_n\setminus E_i$, that is, $x^{\gamma_i}=x^\gamma$ for $x\in E_i$ and $x^{\gamma_i}=x$ for $x\in C_n\setminus E_i$. Thus, $(\a_{s+1})_1\in \Aut(\Gamma)$ and $(\a_{s+1})_1$ fixes $E_0\cup E_2\cup E_3$ pointwise.
% In particular, $(\a_{s+1})_1$ fixes $1\in E_0$.

Suppose that $\Gamma$ is normal. By Proposition~\ref{normal},  $(\a_{s+1})_1\in \Aut(C_n,S)$ as  $(\a_{s+1})_1$ fixes $1$, and since $(\a_{s+1})_1$ fixes $E_0\cup E_2\cup E_3$ pointwise and $|E_0\cup E_2\cup E_3|=3n/4$, it fixes every vertex in $C_n$ as $C_n=\langle E_0\cup E_2\cup E_3\rangle$, which is impossible because $(\a_{s+1})_1$ interchanges the two elements of every coset of $I$ in $E_1$. Thus, $\Gamma$ is non-normal, as required.\qed

\medskip
For a finite group $G$ and a prime $p$, denote by $G_p$ a Sylow $p$-subgroup of $G$. The right regular representation $R(G)$ and the automorphism group $\Aut(G)$ are permutation groups on $G$. Since $R(G)$ is regular on $G$, we have $R(G)\cap \Aut(G)=1$ and $R(G)\Aut(G)=R(G)\rtimes\Aut(G)$, where
\begin{equation}\label{add}
R(g)^\a=R(g^\a) \mbox{ for all } g\in G \mbox{ and } \a\in \Aut(G).
\end{equation}
The normalizer of $R(G)$ in the symmetric group $S_G$ on $G$ is called the {\em holomorph} of $G$, denoted by $\Hol(G)$, and by \cite[Lemma 7.16]{Rot}, $\Hol(G)=R(G)\rtimes\Aut(G)$. Thus, for any $T\leq \Aut(G)$, we have the semiproduct $R(G)\rtimes T$. In particular, $\Hol(C_n)=R(C_n)\rtimes\Aut(C_n)=R(C_n)\rtimes(\Aut(C_{p_1^{k_1}}) \times \Aut(C_{p_2^{k_2}})\times \cdots \times \Aut(C_{p_s^{k_s}})\times \Aut(C_{p_{s+1}^{k_{s+1}}}))$. Note that for all $1\leq i,j\leq s+1$ with $i\not=j$, $\Aut(C_{p_j^{k_j}})$ has a trivial action on $C_{p_i^{k_i}}$. Thus,
\begin{equation}\label{commutes}
[R(C_{p_i^{k_i}}),\Aut(C_{p_j^{k_j}})]=1 \mbox{ and }
[\Aut(C_{p_i^{k_i}}),\Aut(C_{p_j^{k_j}})]=1, \mbox{ for all } i\not=j \end{equation}
that is,  $\Aut(C_{p_j^{k_j}})$ commutes with $R(C_{p_i^{k_i}})$ and $\Aut(C_{p_i^{k_i}})$ pointwise.

Let $\Ga=\Cay(C_n,S)$ be a normal Cayley digraph of $C_n$. Then $\Aut(\Gamma)=R(C_n)\rtimes \Aut(C_n,S)\leq \Hol(C_n)$. Since $\Aut(C_{p_i^{k_i}},S)=\{\a\in \Aut(C_{p_i^{k_i}})\ |\ S^\a=S\}\leq \Aut(C_n,S)$, we have the semiproduct $R(C_n)\rtimes \Aut(C_{p_i^{k_i}},S)\leq \Aut(\Ga)$.

\begin{lemma}\rm\label{RestrictionForH}
Let $\Ga=\Cay(C_n,S)$ be a normal Cayley digraph on the cyclic group $C_n$, and assume that $H$ is a regular subgroup of $\Aut(\Ga)$ isomorphic to $C_n$. Then $H_{p_i}=R(C_{p_i^{k_i}})$ for every $1\leq i\leq s$, and $H\leq R(C_n)\rtimes \Aut(C_{p_{s+1}^{k_{s+1}}},S)$.
\end{lemma}

\proof Let $1\leq i\leq s$. Then $p_i$ is an odd prime. Since $\Aut(C_n)$ is abelian, $\Aut(C_n)_{p_i}$ is the unique Sylow $p_i$-subgroup of $\Aut(C_n)$, and since $R(C_n)\unlhd \Hol(C_n)=R(C_n)\rtimes\Aut(C_n)$, it follows that $\Hol(C_n)_{p_i}\le R(C_n)\Aut(C_n)_{p_i}$. Since $\Aut(C_{p_i^{k_i}})\cong C_{p_i^{k_i-1}(p_i-1)}$,  $\Aut(C_{p_i^{k_i}})$ has a unique subgroup of order $p_i-1$, denoted by  $\Aut(C_{p_i^{k_i}})_{p_i-1}$, and hence  $\Aut(C_{p_i^{k_i}})=\Aut(C_{p_i^{k_i}})_{p_i-1}\times \Aut(C_{p_i^{k_i}})_{p_i}$.
Note that $\Aut(C_{p_{s+1}^{k_{s+1}}})$ is a $2$-group as $p_{s+1}=2$. By Eq~(\ref{factorizationofn}), $p_1>p_2>\cdots>p_s$, and by Eq~(\ref{autoc_nfactorization}), we have
\begin{equation}\label{p_ipart}
\Aut(C_n)_{p_i}\le \Aut(C_{p_1^{k_1}})_{p_1-1}\times \cdots \times \Aut(C_{p_{i-1}^{k_{i-1}}})_{p_{i-1}-1}\times \Aut(C_{p_i^{k_i}})_{p_i}.
\end{equation}

Since $H\cong C_n$, we have $H_{p_i}\cong C_{p_i^{k_i}}$ for every $1\leq i\leq s+1$. Set $A=\Aut(\Ga)$. Since $\Ga=\Cay(C_n,S)$ is normal, Proposition~\ref{normal} implies $A=R(C_n)\rtimes\Aut(C_n,S)\leq
R(C_n)\rtimes \Aut(C_n)=\Hol(C_n)$. For any $h\in H\leq A$, by Eq~(\ref{autoc_nfactorization}) we have

\begin{equation}\label{Helement}
h=R(c)\b_1\cdots\b_s\b_{s+1}, \mbox{ for some } c\in C_n \mbox{ and } \b_i\in \Aut(C_{p_i^{k_i}}).
\end{equation}

Clearly, $R(C_n)\rtimes\Aut(C_{p_{s+1}^{k_{s+1}}},S)\leq R(C_n)\rtimes\Aut(C_{p_{s+1}^{k_{s+1}}})$. Let   $H\leq R(C_n)\rtimes \Aut(C_{p_{s+1}^{k_{s+1}}})$. Since $H\leq A$ and $R(C_n)\leq A$, we have $H\leq A\cap (R(C_n)\rtimes \Aut(C_{p_{s+1}^{k_{s+1}}}))=R(C_n)(A\cap \Aut(C_{p_{s+1}^{k_{s+1}}}))$. Clearly, $A\cap \Aut(C_{p_{s+1}^{k_{s+1}}})=\Aut(C_{p_{s+1}^{k_{s+1}}},S)$. It follows that   $H\leq R(C_n)\rtimes \Aut(C_{p_{s+1}^{k_{s+1}}},S)$. Thus, $H\leq R(C_n)\rtimes \Aut(C_{p_{s+1}^{k_{s+1}}})$ if and only if $H\leq R(C_n)\rtimes \Aut(C_{p_{s+1}^{k_{s+1}}},S)$.  To finish the proof, by Eq~(\ref{Helement}) it suffices to show that $\b_i=1$ and $H_{p_i}=R(C_{p_i^{k_i}})$ for every $1\leq i\leq s$. Let us precess by induction on $i$.

For $i=1$, suppose $H_{p_1}\not=R(C_{p_1^{k_1}})$. Since $R(C_{p_1^{k_1}})$ is characteristic in $R(C_n)$ and $R(C_n)\unlhd A$, we have $R(C_{p_1^{k_1}})\unlhd A$ and hence $R(C_{p_1^{k_1}})H_{p_1}\leq A$. It follows that $R(C_{p_1^{k_1}})H_{p_1}$ is a $p_1$-subgroup of $A$ and $|R(C_{p_1^{k_1}})H_{p_1}|>p_1^{k_1}$. By Proposition~\ref{normal}, $p_1\mid |\Aut(C_n,S)|$ and hence $1\not=\Aut(C_n,S)_{p_1}\leq \Aut(C_n)_{p_1}$. By Eq~(\ref{p_ipart}), $\Aut(C_n)_{p_1}=\Aut(C_{p_1}^{k_1})_{p_1}$, and since $p_1$ is odd, $\Aut(C_{p_1}^{k_1})_{p_1}$ has a unique subgroup of order $p_1$. It follows that $\Aut(C_n)_{p_1}$ and $\Aut(C_n,S)_{p_1}$ contains the unique subgroup of order $p_1$, and by Eq~(\ref{auorderp_i}), we have $\a_1\in \Aut(C_n,S)_{p_1}$, which is impossible by Lemma~\ref{pi}. Therefore, $H_{p_1}=R(C_{p_1^{k_1}})$. This implies that $H\leq C_A(H_{p_1})$ and $R(C_n)\leq C_A(H_{p_1})$, where $C_A(H_{p_1})$ is the centralizer of $H_{p_1}$ in $A$. Since $\b_i\in \Aut(C_{p_i^{k_i}})$, by Eq~(\ref{commutes}) we have $\b_i\in C_A(H_{p_1})$ for every $i\not=1$, and by Eq~(\ref{Helement}) we have $\b_1\in C_A(H_{p_1})$. It follows that $\b_1$ fixes $R(C_{p_1^{k_1}})$ pointwise by conjugacy, and by Eq~(\ref{add}), $\b_1$ fixes $C_{p_1^{k_1}}$ pointwise. By Eq~(\ref{commutes}), $\b_1$ fixes $C_{p_i^{k_i}}$ pointwise for every $2\leq i\leq s+1$, and hence  $\b_1=1$. The inductive basis follows.

Let $2\leq r\leq s$. By inductive hypothesis, we may assume that $\b_i=1$ and $H_{p_i}=R(C_{p_i}^{k_i})$ for every $2\leq i<r$. We only need to show that $\b_r=1$ and $H_{p_r}=R(C_{p_r^{k_r}})$.

From the inductive hypothesis,  Eq~(\ref{Helement}) implies that $H\leq R(C_n)\rtimes (\Aut(C_{p_r^{k_r}})\times \cdots \times \Aut(C_{p_s^{k_s}})\times \Aut(C_{p_{s+1}^{k_{s+1}}}))$, and since $p_2>\cdots>p_s>p_{s+1}=2$, Eq~(\ref{p_ipart}) implies that $(HR(C_n))_{p_r}\leq R(C_n)\rtimes \Aut(C_{p_r^{k_r}})_{p_r}$.

Suppose $H_{p_r}\not=R(C_{p_r^{k_r}})$. Since  $R(C_{p_r^{k_r}})$ is characteristic in $R(C_n)$ and $R(C_n)\unlhd A$, it follows that $R(C_{p_r^{k_r}})\unlhd A$ and $R(C_{p_r^{k_r}})H_{p_r}\leq A$ with $|R(C_{p_r^{k_r}})H_{p_r}|>p_r^{k_r}$.
Then $(HR(C_n))_{p_r}\leq A\cap (R(C_n)\rtimes \Aut(C_{p_r^{k_r}})_{p_r})=R(C_n)(A\cap \Aut(C_{p_r^{k_r}})_{p_r})$. It is easy to see that $A\cap \Aut(C_{p_r^{k_r}})_{p_r}=\Aut(C_{p_r^{k_r}},S)_{p_r}$. Thus,  $R(C_{p_r^{k_r}})H_{p_r}=(HR(C_n))_{p_r}\leq R(C_n)\Aut(C_{p_r^{k_r}},S)_{p_r}$, and since $|R(C_{p_r^{k_r}})H_{p_r}|>p_r^{k_r}$, we have $\Aut(C_{p_r^{k_r}},S)_{p_r}\not=1$.
Since $p_r$ is odd, $\Aut(C_{p_r^{k_r}},S)_{p_r}$ contains the unique subgroup of $\Aut(C_{p_r^{k_r}})$ of order $p_r$. By Eq~(\ref{auorderp_i}), $\a_r\in \Aut(C_{p_r^{k_r}},S)_{p_r}\leq \Aut(C_n,S)$, which is impossible by Lemma~\ref{pi}.

Now we have $H_{p_r}=R(C_{p_r^{k_r}})$, and hence $H\leq C_A(H_{p_r})$ and $R(C_n)\leq C_A(H_{p_r})$. Since $\b_r\in \Aut(C_{p_r^{k_r}})$, by Eq~(\ref{commutes}) we have $\b_r\in C_A(H_{p_i})$ for every $i\not=r$, and by Eq~(\ref{Helement}) we have $\b_r\in C_A(H_{p_r})$. It follows that $\b_r$ fixes $C_{p_r^{k_r}}$ pointwise, and by Eq~(\ref{commutes}), $\b_r$ fixes $C_{p_i^{k_i}}$ pointwise for every $i\not=r$. Thus, $\b_r=1$, as required.
\qed

\medskip

Now we are ready to prove Theorem 1.1.

\medskip

\noindent{\bf Proof of Theorem~\ref{mainth1}:} Let $\Ga=\Cay(C_n,S)$ be a normal Cayley digraph of $C_n$.
Write $A=\Aut(\Ga)$ and $B=R(C_n)\rtimes \Aut(C_{p_{s+1}^{k_{s+1}}},S)$. By Proposition~\ref{normal}, $A=R(C_n)\rtimes\Aut(C_n,S)$. Let $H$ be a regular subgroup of $A$ isomorphic to $R(C_n)$. By Lemma~\ref{RestrictionForH}, $H\leq B$. To finish proof, we only need to show that $H\unlhd A$.

Assume that $\Aut(C_{p_{s+1}^{k_{s+1}}},S)$ has an element $\gamma$ with $o(\gamma)=2^t\geq 8$. Then $k_{s+1}\geq 5$ and $o(\gamma^{2^{t-2}})=4$. Since $\Aut(C_{p_{s+1}^{k_{s+1}}})\cong C_2\times C_{2^{k_{s+1}-2}}$ and $2^{t-2}\geq 2$, we have $\gamma^{2^{t-2}}=\b$ or $\b^{-1}$, where $\b$ is given in Eq~(\ref{auorder4}). This is impossible by Lemma~\ref{pi} because $\gamma^{2^{t-2}}\in\Aut(C_n,S)$.

Now assume that $\Aut(C_{p_{s+1}^{k_{s+1}}},S)$ has no element of order $8$.
Since $\Aut(C_{p_{s+1}^{k_{s+1}}})\cong C_2\times C_{2^{k_{s+1}-2}}$, we  have $\Aut(C_{p_{s+1}^{k_{s+1}}},S)\leq \langle \a,\b\rangle\cong C_2\times C_4$, where  $\a$ is the automorphism of $C_{p_{s+1}^{k_{s+1}}}$ mapping every element to its inverse.  Since $H\cong C_n$ and $H\leq B$, we have $H=\langle R(c)\gamma\rangle$ for some $c\in C_n$ and $\gamma\in \langle \a,\b\rangle$, where $o(\gamma)\leq 4$ and $\gamma\in\Aut(C_{p_{s+1}^{k_{s+1}}},S)$.

Suppose $o(\gamma)=4$. Then $k_{s+1}\geq 4$, and $\Aut(C_{p_{s+1}^{k_{s+1}}})\cong C_2\times C_{2^{k_{s+1}-2}}$ has exactly $4$ elements of order $4$, that is, $\b$, $\b^{-1}$, $\b\a$ and $(\b\a)^{-1}$. By Lemma~\ref{pi}, we may let $\gamma=\b\a$.

Write $E=\langle a_1a_2\cdots a_sa_{s+1}^4\rangle$ and $\Omega=\{E_0,E_1,E_2,E_3\}$, where $E_i=a_{s+1}^iE$ for $i=0,1,2$ or $3$. Then $C_n=E_0\cup E_1\cup E_2\cup E_3$. Since $R(C_n)$ is cyclic, $R(E)$ is characteristic in  $R(C_n)$, and since $R(C_n)\unlhd A$, we have $R(E)\unlhd A$. Clearly, $\Omega$ is the orbit set of $R(E)$ on $V(\Gamma)$, and hence $\Omega$ is a complete imprimitive block system of $A$. Note that $E_0\cup E_2=\langle a_1a_2\cdots a_sa_{s+1}^2\rangle$ is also a characteristic subgroup of $C_n$. Similarly, $\Omega':=\{E_0\cup E_2,E_1\cup E_3\}$ is a complete imprimitive block system of $A$. Since $H$ is regular on $V(\Gamma)=C_n$, $H$ is transitive on both $\Omega$ and $\Omega'$. Note that $\b$ fixes every $E_i$ setwise, and $\a$ fixes $E_0$ and $E_2$, but interchanges $E_1$ and $E_3$. It follows that $\g$ fixes $E_0$ and $E_2$, but interchanges $E_1$ and $E_3$, implying that $\g$ induces a trivial permutation on $\Omega'$.

Recall $H=\langle R(c)\gamma\rangle$. If $R(c)$ induces a trivial action on $\Omega'$, then $H$ induces a trivial action on $\Omega'$ because  $\gamma$ induces a trivial action on $\Omega'$, contradicting that $H$ is transitive on $\Omega'$. Thus, $R(c)$ is transitive on $\Omega'$, implying that $R(c)$ maps $E_0$ to $E_1$ or $E_3$. It follows that $R(c)=R(a)^s$ for some odd integer $s$. Then the induced permutation of $R(c)$ on $\Omega$ is $(E_0\ E_1\ E_2\ E_3)$ or $(E_0\ E_3\ E_2\ E_1)$.  Since $\g$ induces the permutation $(E_1\ E_3)$ on $\Omega$, $R(c)\g$ induces an permutation of order $2$ on $\Omega$, contradicting the transitivity of $H$ on $\Omega$.

Thus, $o(\gamma)\not=4$. It follows that $\gamma^2=1$ and $(R(c)\gamma)^2\in R(C_n)\gamma R(C_n)\gamma=R(C_n)$ as $R(C_n)\unlhd A$.  Clearly, $\langle (R(c)\gamma)^2\rangle$ is a subgroup of order $n/2$ in $H$, and hence $\langle (R(c)\gamma)^2\rangle=\langle R(a^2)\rangle$, the unique subgroup of order $n/2$ in $R(C_n)$. Write $D=\langle R(a^2)\rangle$. Since $\langle R(a^2)\rangle$ is characteristic in $R(C_n)$, $D$ is normal in $A$ and $D\leq H\cap R(C_n)$. Recall that $A=R(C_n)\rtimes \Aut(C_n,S)$. Then $R(C_n)/D$ is a normal subgroup of order $2$ in $A/D$, and therefore $R(C_n)/D$ lies in the center of $A/D$. Thus, $A/D=R(C_n)/D\times \Aut(C_n,S)D/D$, and since $\Aut(C_n,S)$ is abelian, every subgroup of $A/D$ is normal. In particular, $H/D\unlhd A/D$, implying $H\unlhd A$, as required. \qed

\section{Proof of Theorem~\ref{mainth2}\label{theorem2}}

Let $\D_{2n}$ be the dihedral group of order $2n$. Write
\begin{equation}\label{d_2n}
\D_{2n}=\lg a,b\ |\ a^n=b^2=1,b^{-1}ab=a^{-1}\rg\ \mbox{ and }
C=\langle a\rangle.
\end{equation}
    Then $C$ is a cyclic subgroup of order $n$, which is characteristic in $\D_{2n}$ when $n\geq 3$.

By Proposition~\ref{NDCID}, if $n$ is odd or $n=2$ or $4$, then $\D_{2n}$ is an NDCI-group and NCI-group, which implies that $\D_{2n}$ is neither an NNN-group nor an NNND-group. Thus, we may assume that $n\geq 6$ is even, and for convenience to prove Theorem~\ref{mainth2}, it is divided into two parts: 1. $n$ is even and $n\geq 6$ with $n\not=8$; 2. $n=8$. For the first part, we have the following result.

\begin{lemma}\rm\label{mainlem} Let $n$ be even and let $n\geq 6$ with $n\not=8$.
Set $S=\{a,a^{-1},b,a^{\frac{n}{2}}b\}$ when $n/2$ is odd, and $S=\{a,a^{-1},b,a^{\frac{n}{4}}b,a^{\frac{n}{2}}b,a^{\frac{3}{4}n}b\}$ when $n/2$ is even. Then $\Ga=\Cay(\D_{2n},S)$ is an NNN graph.
\end{lemma}

\proof  It is easy to see that $S=S^{-1}$ and $\langle S\rangle=\D_{2n}$. Then $\Ga$ is a connected graph. For a non-negative integer $i$ and $u\in V(\Ga)$, denote by  $\Ga_i(u)$  the $i$-th neighborhood of $u$ in $\Ga$, that is, $\Ga_i(u)=\{v\in V(\Ga)\ |\ d(u,v)=i\}$, where $d(u,v)$ is the distance between $u$ and $v$ in $\Ga$. Furthermore, write  $\Ga_i[u]=\{v\in V(\Ga)\ |\ d(u,v)\leq i\}$, the set of vertices with distance at most $i$ from $u$. For short, write $\Ga(u)=\Ga_1(u)$. Let $A_1$ be the stabilizer of 1 in $A$ and let $A_1^{\ast}$ be the subgroup of $A_1$ fixing the neighbourhood  $\Ga(1)$ of 1 in $\Ga$ pointwise.

By Eq~(\ref{d_2n}), it is easy to see that $a^{-1}$ and $b$ has the same relations as $a$ and $b$ in $\D_{2n}$. Since $\langle a^{-1},b\rangle=\D_{2n}$, we denote by $\a$ be the automorphism of $\D_{2n}$ induced by
$$\a: a\mapsto a^{-1},\ \ b\mapsto b.$$
Similarly, let $\b$ and $\g$ (if $n$ has a divisor $4$) be the automorphisms of $\D_{2n}$ induced by
$$\b: a\mapsto a,\ b\mapsto a^{n/2}b;\ \ \ \ \g: a\mapsto a, \ b\mapsto a^{n/4}b \ (4 \mbox{ divides } n). $$

It is easy to see that $\langle \a,\b\rangle\cong C_2\times C_2$ and $\langle \a,\g\rangle\cong \D_8$ with $\g^\a=\g^{-1}$. Write $A=\Aut(\Ga)$. Now we determine $A$ by considering the parity of $n/2$.

\medskip
\noindent{\bf Claim 1:} If $n/2$ is odd, then $A=R(\D_{2n})\rtimes \Aut(\D_{2n},S)=R(\D_{2n})\rtimes \langle \a,\b\rangle$.

Let $n/2$ be odd. Then $S=\{a,a^{-1},b,a^{\frac{n}{2}}b\}$, and it is easy to see that $\a,\b\in \Aut(\D_{2n},S)$. Thus, $|\Aut(\D_{2n},S)|\geq |\langle \a,\b\rangle|=|C_2\times C_2|=4$. Note that $A=R(\D_{2n})A_1$ and $\Aut(\D_{2n},S)\leq A_1$. To prove Claim~1, it suffices to show that $|A_1|\leq 4$.

\begin{figure}[htb]
\begin{center}
\begin{tikzpicture}[node distance=1.2cm,thick,scale=0.7,every node/.style={transform shape},scale=1.2](330,180)(20,30)
%%\thicklines

\draw(0,10)--(-3,8);\draw(0,10)--(3,8);\draw(0,10)--(-7,8);\draw(0,10)--(7,8);
\draw(-3,8)--(-5,4);\draw(-3,8)--(0,4);\draw(-3,8)--(3,4);
\draw(-7,8)--(-7,4);\draw(-7,8)--(-5,4);\draw(-7,8)--(-3,4);
\draw(3,8)--(5,4);\draw(3,8)--(0,4);\draw(3,8)--(-3,4);
\draw(7,8)--(7,4);\draw(7,8)--(5,4);\draw(7,8)--(3,4);

\draw(0,10)node [above]{1};
\draw(-3,8)node [above]{$b$};
\draw(3,8)node [above]{$a^{\frac{n}{2}}b$};
\draw(-7,8)node [above]{$a$};
\draw(7,8)node [above]{$a^{-1}$};
\draw(0,4)node [below]{$a^{\frac{n}{2}}$};
\draw(-7,4)node [below]{$a^2$};
\draw(7,4)node [below]{$a^{-2}$};
\draw(-5,4)node [below]{$ba$};
\draw(5,4)node [below]{$a^{\frac{n}{2}+1}b$};
\draw(-3,4)node [below]{$a^{\frac{n}{2}-1}b$};
\draw(3,4)node [below]{$ab$};

\filldraw [black] (0,10) circle (2pt);
\filldraw [black] (-3,8) circle (2pt);
\filldraw [black] (-7,8) circle (2pt);
\filldraw [black] (3,8) circle (2pt);
\filldraw [black] (7,8) circle (2pt);
\filldraw [black] (0,4) circle (2pt);
\filldraw [black] (-7,4) circle (2pt);
\filldraw [black] (7,4) circle (2pt);
\filldraw [black] (-5,4) circle (2pt);
\filldraw [black] (5,4) circle (2pt);
\filldraw [black] (-3,4) circle (2pt);
\filldraw [black] (3,4) circle (2pt);
\end{tikzpicture}
\end{center}
\caption{The induced subgraph $[\Ga_2[1]]$}
\label{figure1}
\end{figure}
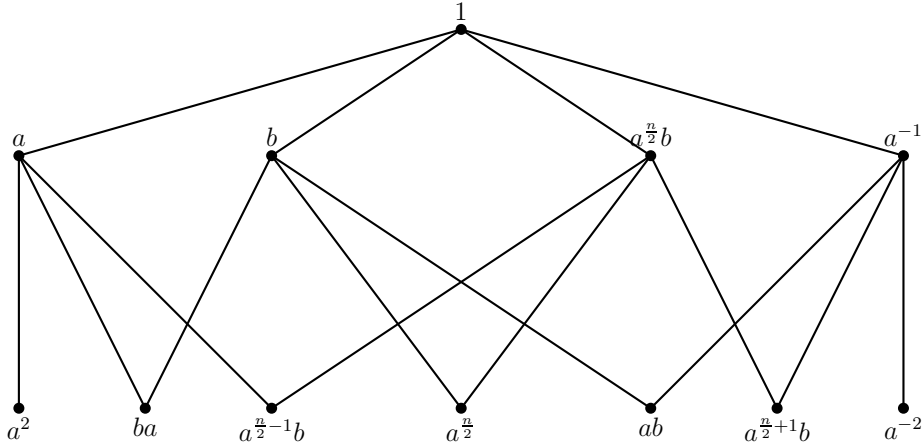

Since $n\ge 6$, the induced subgraph  $[\Ga_2[1]]$ of  $\Ga_2[1]$ in $\Ga$ can be drawn as Figure~\ref{figure1}.

 Note that $A_1^{\ast}$ fixes $\Ga_i[1]$ pointwise for $i=0$ and $1$, and fixes $\Ga_2(1)$ setwise. By Figure~\ref{figure1}, all vertices on $[\Ga_2[1]]$  are distinct, and $\Ga(x)\cap \Ga(1)$ for all $x\in \Ga_2(1)$ are different from one another, implying that  $A_1^{\ast}$ fixes $\Ga_2(1)$ pointwise. By the connectivity of $\Ga$,  $A_1^{\ast}$ fixes every vertex in $\Ga$, that is, $A_1^*=1$. Thus, $A_1$ has a faithful action on $S$.

 By Figure~\ref{figure1}, for every vertex $x\in \Ga(b)\cap \Ga_2(1)$ there exist a $4$-cycle passing through $1$, $b$ and $x$, and for every vertex $y\in \Ga(a^{n/2}b)\cap \Ga_2(1)$, there exist a $4$-cycle passing through $1$, $a^{n/2}b$ and $y$. However, there is no $4$-cycle passing through either $1$, $a$ and $a^2\in \Ga(a)\cap \Ga_2(1)$, or $1$, $a^{-1}$ and $a^{-2}\in \Ga(a^{-1})\cap \Ga_2(1)$. Since $A_1$ fixes $\Ga_i(1)$ for each non-negative integer $i$, $A_1$ fixes $\{a,a^{-1}\}$ and $\{b,a^{n/2}b\}$ setwise, respectively. Since $A_1$ acts faithfully on $S$, we have that $|A_1|\leq |\{a,a^{-1}\}|!\cdot |\{b,a^{n/2}b\}|!=4$, as required.
\qed
\medskip

\medskip
\noindent{\bf Claim 2:} If $n/2$ is even, then $A=R(\D_{2n})\rtimes \Aut(\D_{2n},S)=R(\D_{2n})\rtimes\langle \a,\g\rangle$.

Let $n/2$ be even. Then $S=\{a,a^{-1},b,a^{\frac{n}{4}}b,a^{\frac{n}{2}}b,a^{\frac{3}{4}n}b\}$, and it is easy to see that $\a,\g\in \Aut(\D_{2n},S)$. Thus, $|\Aut(\D_{2n},S)|\geq |\langle \a,\g\rangle|=|\D_8|=8$. Since $A=R(\D_{2n})A_1$ and $\Aut(\D_{2n},S)\leq A_1$, to prove Claim~2 we only need to show that $|A_1|\leq 8$.

Write $S_1=\{a,a^{-1}\}$ and $S_2=\{b,a^{\frac{n}{4}}b,a^{\frac{n}{2}}b,a^{\frac{3}{4}n}b\}$. Then $S=S_1\cup S_2$ and $\Ga=\Ga_1\cup \Ga_2$, where $\Ga_1=\Cay(\D_{2n},S_1)$ and  $\Ga_2=\Cay(\D_{2n},S_2)$. Since $n/2$ is even and $n\geq 6$ with $n\not=8$, we have $n\geq 12$. For convenience, we say that an edge in $\Ga_1$ is an {\em $a$-edge} of $\Ga$, and an edge in $\Ga_2$ is a {\em $b$-edge} of $\Ga$. Clearly, $\Ga_1$ is a union of two $n$-cycles, that is, the induced subgraph $[C]$ and $[Cb]$ of $C$ and $Cb$ in $\Ga$, where $\D_{2n}=C\cup Cb$.

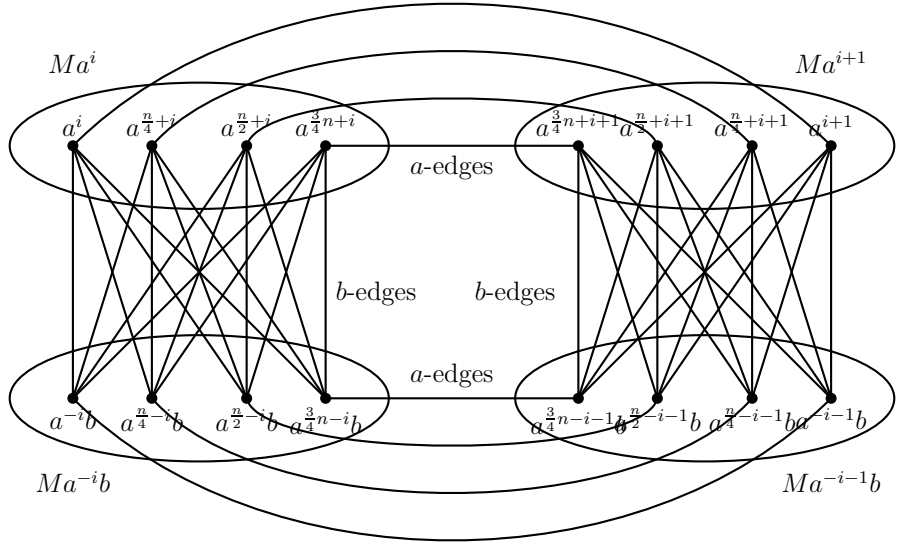
\begin{figure}[htb]
\begin{center}
\begin{tikzpicture}[node distance=1.2cm,thick,scale=0.7,every node/.style={transform shape},scale=1.2](330,180)(20,30)
%\tikzstyle{every node}=[font=\tiny]
\thicklines

\draw(-4,0) ellipse [x radius=3,y radius=1];
\draw(4,0) ellipse [x radius=3,y radius=1];
\draw(-4,-4) ellipse [x radius=3,y radius=1];
\draw(4,-4) ellipse [x radius=3,y radius=1];

\draw(-6,0)--(-6,-4);\draw(-6,0)--(-4.75,-4);\draw(-6,0)--(-3.25,-4);\draw(-6,0)--(-2,-4);
\draw(-4.75,0)--(-6,-4);\draw(-4.75,0)--(-4.75,-4);\draw(-4.75,0)--(-3.25,-4);\draw(-4.75,0)--(-2,-4);
\draw(-3.25,0)--(-6,-4);\draw(-3.25,0)--(-4.75,-4);\draw(-3.25,0)--(-3.25,-4);\draw(-3.25,0)--(-2,-4);
\draw(-2,0)--(-6,-4);\draw(-2,0)--(-4.75,-4);\draw(-2,0)--(-3.25,-4);\draw(-2,0)--(-2,-4);

\draw(6,0)--(6,-4);\draw(6,0)--(4.75,-4);\draw(6,0)--(3.25,-4);\draw(6,0)--(2,-4);
\draw(4.75,0)--(6,-4);\draw(4.75,0)--(4.75,-4);\draw(4.75,0)--(3.25,-4);\draw(4.75,0)--(2,-4);
\draw(3.25,0)--(6,-4);\draw(3.25,0)--(4.75,-4);\draw(3.25,0)--(3.25,-4);\draw(3.25,0)--(2,-4);
\draw(2,0)--(6,-4);\draw(2,0)--(4.75,-4);\draw(2,0)--(3.25,-4);\draw(2,0)--(2,-4);

\draw(-6,0)..controls (-3,3) and (3,3)..(6,0);
\draw(-4.75,0)..controls (-3,2) and (3,2)..(4.75,0);
\draw(-3.25,0)..controls (-3,1) and (3,1)..(3.25,0);
\draw(-2,0)--(2,0);
\draw(-6,-4)..controls (-3,-7) and (3,-7)..(6,-4);
\draw(-4.75,-4)..controls (-3,-6) and (3,-6)..(4.75,-4);
\draw(-3.25,-4)..controls (-3,-5) and (3,-5)..(3.25,-4);
\draw(-2,-4)--(2,-4);

\draw(-2,0)node [above]{$a^{\frac{3}{4}n+i}$};
\draw(-3.25,0)node [above]{$a^{\frac{n}{2}+i}$};
\draw(-4.75,0)node [above]{$a^{\frac{n}{4}+i}$};
\draw(-6,0)node [above]{$a^i$};
\draw(2,0)node [above]{$a^{\frac{3}{4}n+i+1}$};
\draw(3.25,0)node [above]{$a^{\frac{n}{2}+i+1}$};
\draw(4.75,0)node [above]{$a^{\frac{n}{4}+i+1}$};
\draw(6,0)node [above]{$a^{i+1}$};
\draw(-2,-4)node [below]{$a^{\frac{3}{4}n-i}b$};
\draw(-3.25,-4)node [below]{$a^{\frac{n}{2}-i}b$};
\draw(-4.75,-4)node [below]{$a^{\frac{n}{4}-i}b$};
\draw(-6,-4)node [below]{$a^{-i}b$};
\draw(2,-4)node [below]{$a^{\frac{3}{4}n-i-1}b$};
\draw(3.25,-4)node [below]{$a^{\frac{n}{2}-i-1}b$};
\draw(4.75,-4)node [below]{$a^{\frac{n}{4}-i-1}b$};
\draw(6,-4)node [below]{$a^{-i-1}b$};
\draw(-6,1)node [above]{$Ma^i$};
\draw(6,1)node [above]{$Ma^{i+1}$};
\draw(-6,-5)node [below]{$Ma^{-i}b$};
\draw(6,-5)node [below]{$Ma^{-i-1}b$};

\draw(0,-3.3)node [below]{$a$-edges};
\draw(0,0)node [below]{$a$-edges};
\draw(-1.2,-2)node [below]{$b$-edges};
\draw(1,-2)node [below]{$b$-edges};

\filldraw [black] (-2,0) circle (2pt);
\filldraw [black] (-3.25,0) circle (2pt);
\filldraw [black] (-4.75,0) circle (2pt);
\filldraw [black] (-6,0) circle (2pt);
\filldraw [black] (2,0) circle (2pt);
\filldraw [black] (3.25,0) circle (2pt);
\filldraw [black] (4.75,0) circle (2pt);
\filldraw [black] (6,0) circle (2pt);
\filldraw [black] (-2,-4) circle (2pt);
\filldraw [black] (-3.25,-4) circle (2pt);
\filldraw [black] (-4.75,-4) circle (2pt);
\filldraw [black] (-6,-4) circle (2pt);
\filldraw [black] (2,-4) circle (2pt);
\filldraw [black] (3.25,-4) circle (2pt);
\filldraw [black] (4.75,-4) circle (2pt);
\filldraw [black] (6,-4) circle (2pt);
\end{tikzpicture}
\end{center}
\caption{The induced subgraph $[Ma^i\cup Ma^{i+1}\cup Ma^{-i}b\cup  Ma^{-i-1}b]$.}
\label{figure2}
\end{figure}

Let $M=\langle a^{\frac{n}{4}}\rangle$. Then $M\leq C$ is a subgroup of order $4$, and hence characteristic in $\D_{2n}$.
Note that $\D_{2n}=\cup_{i=0}^{n/4-1}Ma^i\cup Ma^ib$ and the induced subgraphs $[Ma^i]$ and $[Ma^ib]$ in $\Ga$ have no $a$-edges, and also have no $b$-edges because $S_2=Mb$. Thus, $[Ma^i]$ and $[Ma^ib]$ are empty graphs of order $4$.  Since $n\geq 12$, the induced subgraphs $[Ma^i\cup Ma^{i+1}]$ and $[Ma^{-i}b\cup Ma^{-i-1}b]$ are matchings, and  the induced subgraphs $[Ma^i\cup Ma^{-i}b]$ and  $[Ma^{i+1}\cup Ma^{-i-1}b]$
are isomorphic to $K_{4,4}$, where the powers of $a$ are taken modulo $n/4$. Thus, the induced subgroup $[Ma^i\cup Ma^{i+1}\cup Ma^{-i}b\cup  Ma^{-i-1}b]$ can be drawn as Figure~\ref{figure2}.

It is easy to see that $\Ga=\cup_{i=0}^{n/4-1}[Ma^i\cup Ma^{i+1}\cup Ma^{-i}b\cup  Ma^{-i-1}b]$. By Figure~\ref{figure2}, there is exactly four $4$-cycle passing through every $a$-edge, and there are at least nine $4$-cycles passing through every $b$-edge. This implies that $A$ fixes the set of $a$-edges of $\Ga$ setwise, and the set of $b$-edges of $\Ga$ setwise, forcing that $A\leq \Aut(\Ga_1)$ and $A\leq \Aut(\Ga_2)$.

Since $A\leq \Aut(\Ga_1)$, $A_1$ fixes the $n$-cycle $[C]$ and hence $|A_1|=|a^{A_1}|\cdot|A_{1a}|\leq 2|A_{1a}|$, where $A_{1a}$ is the subgroup of $A_1$ fixing $a$ and $a^{A_1}\subseteq \{a,a^{-1}\}$. Since $[C]$ is a cycle, $A_{1a}$ fixes $C$ pointwise. In particular,  $A_{1a}$ fixes $Ma^i$ pointwise for every $0\leq i\leq n/4-1$. By Figure~\ref{figure2}, $A_{1a}$ fixes $Ma^ib$ setwise for every $0\leq i\leq n/4-1$ because  $A_{1a}\leq \Aut(\Ga_2)$ fixes the set of $b$-edges. In particular, $A_{1a}$ fixes $Mb$ setwise, and hence $|A_{1a}|=|b^{A_{1a}}||A_{1ab}|\leq |Mb||A_{1ab}|\leq 4|A_{1ab}|$, where $A_{1ab}$ is the subgroup of $A_{1a}$ fixing $b$. Since $A_{1ab}$ fixes $Ma^{-1}b$ setwise and the induced subgraph $[Mb\cup Ma^{-1}b]$ is a matching, it fixes $a^{-1}b$. Then $A_{1ab}$ fixes the arc $(b,a^{-1}b)$, and since $A_{1ab}\leq \Aut(\Ga_1)$ fixes the $n$-cycle $[Cb]$, $A_{1ab}$ fixes $Cb$ pointwise. It follows that $A_{1ab}=1$ as $A_{1ab}\leq A_{1a}$ fixes $C$ pointwis. Therefore, $A_1\leq 2|A_{1a}|\leq 8|A_{1ab}|=8$, as required.
\qed
\medskip

By Claims~1 and 2, $\a\in A$ and Proposition~\ref{normal} implies that $\Ga$ is normal.

\medskip
\noindent {\bf Claim 3:} $\lg R(ab)\a,R(b)\rg$ is a regular subgroup of $A$ isomorphic to $\D_{2n}$.

Since $a^\a=a^{-1}$ and $b^\a=b$, we have $R(b)^{\a}=R(b)$ and $R(ab)^{\a}=R((ab)^{\a})=R(a^{-1}b)$. Furthermore, $\a$ has order $2$ and
\begin{equation*}
  (R(ab)\a)^{R(b)}=R(ab)^{R(b)}\a=R(b^{-1}abb)\a=R(b^{-1}a)\a=\a R(ba^{-1})=(R(ab)\a)^{-1}.
\end{equation*}
Since $(R(ab)\a)^2=(R(ab)\a)(R(ab)\a)=R(ab)R(a^{-1}b)=R(a^2)$, $o(R(ab)\a)$ is even and
$n/2=o(R(a^2))=o(R(ab)\a)/2$, that is, $o(R(ab)\a)=n$.
Thus, $\lg R(ab)\a,R(b)\rg\cong \D_{2n}$.

To prove the regularity of $\lg R(ab)\a,R(b)\rg$, it suffices to show that it is transitive on $\D_{2n}$.
Since $\lg R(a^2)\rg$ is characteristic in $R(\D_{2n})$ and $R(D_{2n})\unlhd A$, it follows that $\lg R(a^2)\rg\unlhd A$.
Note that $\lg R(a^2)\rg$ has four orbits on $\D_{2n}$: $\lg a^2\rg$, $\lg a^2\rg a$, $\lg a^2\rg b$ and $\lg a^2\rg ab$.
Clearly, $R(b)$ interchanges $\lg a^2\rg$ and $\lg a^2\rg b$, and $\lg a^2\rg a$ and $\lg a^2\rg ab$, respectively. Moreover,  $R(ab)\a$ interchanges $\lg a^2\rg$ and $\lg a^2\rg ab$, and $\lg a^2\rg a$ and $\lg a^2\rg b$, respectively.
It follows that $\lg R(ab)\a,R(b)\rg$ is transitive on $\D_{2n}$, as required.
\qed
\medskip

Let $H=\lg R(ab)\a,R(b)\rg$. To finish the proof, by Claim~3 we only need to show that $H$ is not normal in $A$. Suppose to the contrary that $H\unlhd A$.

Assume $n/2$ is odd. By Claim~1, $\b\in A$, and then
$R(b)^\b R(b)=R(a^{\frac{n}{2}}b)R(b)=R(a^{\frac{n}{2}})\in H$.
Recall that $R(a^2)=(R(ab)\a)^2\in H$. Since $\gcd(\frac{n}{2},2)=1$, we have $R(a)\in H$, and hence $H=\lg R(a),R(ab)\a,R(b)\rg$, forcing $\a\in H$, contradicting the regularity of $H$.

Assume $n/2$ is even. By Claim~2, $\g\in A$, and then
$(R(a)\a)^{\g}=R(a)\a^{\g}=R(a)\a\g^2\in H$. Since $R(a)\a\in H$, we have $\g^2\in H$, which contradicts the regularity of $H$ because $\g^2\neq 1$ fixes $1$. This completes the proof.
\qed
\medskip

Now we deal with the second part: $n=8$.

\begin{lemma}\rm\label{lem2}
$\D_{16}$ is not an NNND-group.
\end{lemma}

\proof Let $\Ga=\Cay(\D_{16},S)$ be a normal Cayley digraph and $A=\Aut(\Ga)$. By Proposition~\ref{normal}, $A=R(\D_{16})\rtimes \Aut(\D_{16},S)$.
Let $H$ be a regular subgroup of $A$ isomorphic to $\D_{16}$. To finish the proof, we only need to show that $H\unlhd A$. If $|\Aut(\D_{16},S)|\leq 2$ then $|A:H|\leq 2$, implying $H\unlhd A$. Thus, we assume that $|\Aut(\D_{16},S)|\geq 4$.

Recall that $\D_{16}=\lg a,b\ |\ a^8=b^2=1,b^{-1}ab=a^{-1}\rg$ and $C=\langle a\rangle$. Let $\a$, $\b$ and $\g$ be automorphisms of $\D_{16}$ induced by
$$\a: a\mapsto a^{-1},\ b\mapsto b;\ \ \b: a\mapsto a^5,\ b\mapsto b;\ \ \g: a\mapsto a, \ b\mapsto ab.$$
By \cite[Lemma 3.1]{XFZ}, $\Aut(\D_{16})=\lg \g\rg\rtimes \lg\a,\b\rg\cong C_8\rtimes (C_2\times C_2)$ with $o(\a)=o(\b)=2$ and $o(\g)=8$, where
\begin{equation}\label{AutD_16}
 \a\b=\b\a,\ \ \g^\a=\g^{-1},\ \ \g^\b=\g^5.
\end{equation}
Clearly, $\Aut(\D_{16})$ has no element of order $16$. The subgroup $\lg\g,\a\rg$ is the dihedral group of order $16$, and $\lg\g^2,\g\a\b\rg$ is the quaternion group of order $8$.

\medskip
\noindent {\bf Observation:}  Let $O_i$ be the set of elements of order $i$ in $\Aut(\D_{16})$. Then $O_2=\lg\g\rg\a\cup \lg\g^2\rg\a\b\cup \{\b,\g^4,\g^4\b\}$, $O_4=\{\g^2, \g^6, \g^2\b, \g^6\b, \g\a\b, \g^3\a\b, \g^5\a\b, \g^7\a\b\}$ and $O_8=\lg \g^2\rg \g\cup \lg \g^2\rg \g\b$. Furthermore, an elementary abelian subgroup of $\Aut(\D_{16})$ containing $\g^i\a$ with $i$ odd, must be $\langle \g^i\a\rg$ or $\langle\g^i\a,\g^4\rg=\{1,\g^i\a,\g^{i+4}\a,\g^4\}$.

\medskip
The first part of Observation follows from the relations in Eq~(\ref{AutD_16}), and the second part follows from the fact that an elementary abelian subgroup in the dihedral subgroup $\lg\g,\a\rg$ containing $\g^i\a$ with odd $i$, must be $\langle \g^i\a\rg$ or $\langle\g^i\a,\g^4\rg$, and $\g^i\a$ cannot commutes with every involution in $O_2\cap \lg\g,\a\rg\b=\lg\g^2\rg\a\b\cup \{\b,\g^4\b\}$.

\medskip
The normalizer of $R(\D_{16})$ in the symmetric group on $\D_{16}$ is called the {\em holomorph} of $\D_{16}$, denoted by $\Hol(\D_{16})$, and by \cite[Lemma 7.16]{Rot}, $\Hol(\D_{16})=R(\D_{16})\rtimes \Aut(\D_{16})$. It follows that  $A\leq \Hol(\D_{16})$. Note that
$R(d)^\d=R(d^\d)$ for all $d\in \D_{16}$ and $\d\in \Aut(\D_{16})$. This, together with Eq~(\ref{AutD_16}), implies that $R(a)\g=\g R(a), R(a)^{R(b)}=R(a)^{-1}, \g^{R(b)}=\g R(a), R(a)^{\a}=R(a)^{-1}, \g^{\a}=\g^{-1}, R(a)^{\b}=R(a)^5, \g^{\b}=\g^5, \a R(b)= R(b)\a, \b R(b)=R(b)\b$ and $\a\b=\b\a$. Then we may set
$$A1=\langle R(a),\g\rg\cong C_8^2,\ A2=\langle \a,\b,R(b)\rg\cong C_2^3, \ B=\langle \a,\b\rg\cong C_2^2 \mbox{ and } D=A1\rtimes B.$$ Furthermore,
\begin{equation}\label{HolD_16C}
\Hol(\D_{16})=D\rtimes\langle R(b)\rangle=A1\rtimes A2\cong C_8^2\rtimes C_2^3.
\end{equation}

For $\d\in B$, we have $R(a)^\d=R(a)^r$ for some odd integer $r$, and
\begin{equation}\label{r-B}
r=1,3,5,7\ (\mod 8) \mbox{ if and only if } \d=1,\a\b,\b,\a, \mbox{ respectively.}
\end{equation}

Recall that $R(a)^{\a}=R(a)^{-1}, \g^{\a}=\g^{-1}, R(a)^{\b}=R(a)^5$, $\g^{\b}=\g^5$, $A1=\lg R(a),\g\rg$ and $B=\lg\a,\b\rg$. This implies that
\begin{equation}\label{ab}
\mbox{if } R(a)^\d=R(a)^r \mbox{ for } \d\in B, \mbox{ then } z^{\d}=z^r \mbox{ for all } z\in A1.
\end{equation}

An element $z\in\Hol(\D_{16})$ is called {\em semiregular}, if all orbits of $z$ on $\D_{16}$ have the same length. Clearly, any power of a semiregular element is semiregular. Since $R(\D_{16})$ and $H$ are regular, every element in these two groups is semiregular. Now we consider some semiregular elements in $\Hol(\D_{16})$.
In what follows, all equations on integers are taken $\mod 8$, and we omit it for convenience. For example, write $r=1,3,5,7\ (\mod 8)$ as $r=1,3,5,7$.
Note that every subgroup of $R(C)$ is characteristic in $R(\D_{16})$ and so normal in $\Hol(\D_{16})$ as $R(\D_{16})\unlhd \Hol(\D_{16})$, implying that its orbit set is a complete imprimitive block system of $\Hol(\D_{16})$. In particular, $\{C,Cb\}$ and $\{\lg a^2\rg,\lg a^2\rg a,\lg a^2\rg b,\lg a^2\rg ab\}$ are  complete imprimitive block systems of $\Hol(\D_{16})$, respectively. By Eq~(\ref{HolD_16C}),  $\Hol(\D_{16})=D\cup DR(b)$.

\medskip
\noindent {\bf Claim 1:} Let $\d\in B$ and $R(a)^\d=R(a)^r$ for some odd integer $r$. Then an element $R(a)^i\g^j\d\in D$ is a semiregular element of order $2$ interchanging $\lg a^2\rg$ and $\lg a^2\rg a$ if and only if
$i$ is odd and $j$ is even with $\d=\a$, and is a semiregular element of order $8$ if and only if  $i$ is odd and $j$ is even with $\d=1$ or $\b$. An element $R(a)^i\g^j\d R(b)\in DR(b)$ is a semiregular element of order $2$ if and only if $(1+r)j=0$ and $j=(1-r)i$, and is a semiregular element of order $8$ if and only if $j$ is even and $(i-ri+rj)/2$ is odd.

\medskip
Let $x=R(a)^i\g^j\d\in D$ with $\d\in B$, and $R(a)^\d=R(a)^r$. Then $r=1,3,5$ or $7$, and by Eq~(\ref{ab}), we have
\begin{equation} \label{square}
x^2=R(a)^i\g^j(R(a)^i\g^j)^{\d}=R(a)^{i+ri}\g^{j+jr}.
\end{equation}

Assume that $x$ is semiregular of order $2$ interchanging $\lg a^2\rg$ and $\lg a^2\rg a$. Then $(1+r)i=0$ and $(1+r)j=0$, and since $\g^j\d$ fixes $\lg a^2\rg$ and $\lg a^2\rg a$ respectively, $R(a)^i$ interchanges $\lg a^2\rg$ and $\lg a^2\rg a$, forcing that $i$ is odd. It follows from $(1+r)i=0$ that $r=7$ and hence $\d=\a$ by Eq~(\ref{r-B}). If $j$ is odd, then $i-j$ is even, say $i-j=2t$, and then $x$ fixes $a^tb$, contradicting the semiregularity of $x$. Thus, $j$ is even.

On the other hand, $R(a)^i\g^j\a$ with $i$ odd and $j$ even, has order $2$ as $(R(a)^i\g^j)^\a=(R(a)^i\g^j)^{-1}$, and interchanges  $\lg a^2\rg$ and $\lg a^2\rg a$. Since $i-j$ is odd,  it is easy to see that $R(a)^i\g^j\a$ has no fixed point on $Cb$, and so is semiregular.

Assume that $x$ is semiregular of order $8$. Then $x^2$ has order $4$. If $o(R(a)^{i+ri})=1$ or $2$, then Eq~(\ref{square}) implies that $x^4=(\g^{j+jr})^2$ has order $2$ that fixes $1$, which is impossible as $x^4$ is semiregular. Thus,  $o(R(a)^{i+ri})=4$, forcing that $i$ is odd and $r=1$ or $5$, and by Eq~(\ref{r-B}), $\d=1$ or $\b$. If $j$ is odd, then $x^4=R(a)^4\g^4$ that fixes $Cb$ pointwise, contradicting the semiregularity of $x^4$. Thus, $j$ is even.

On the other hand, assume that $i$ is odd and $j$ is even with $\d=1$ or $\b$.
By Eq~(\ref{r-B}), $r=1$ or $5$, and hence $o(R(a)^{i+ri})=4$  and $o(x)=8$ by Eq~(\ref{square}), because $o(\g^{j+jr})=1$ or $2$. Furthermore, $x^4=R(a)^{2(i+ri)}=R(a)^4$ that has order $2$, and all orbits of $\lg x^4\rg$ on $\D_{16}$ have length $2$, implying $x$ is semiregular.

Now let $y=R(a)^i\g^j\d R(b)$ with $\d\in B$, and $R(a)^\d=R(a)^r$. Then $r=1,3,5$ or $7$, and by Eq~(\ref{ab}) we have
\begin{equation}
\label{square2}
y^2=R(a)^i\g^j(R(a)^i\g^j)^{R(b)\d}=R(a)^i\g^j(R(a)^{j-i}\g^j)^\d
=R(a)^{i-ri+rj}\g^{j+jr}.
\end{equation}

Assume that $y$ is semiregular of order $2$. Then $(1+r)j=0$ and $i-ri+rj=0$, that is, $(1+r)j=0$ and $j=(1-r)i$. On the other hand, assume that  $(1+r)j=0$ and $i-ri+rj=0$. By Eq~(\ref{square2}), $y$ has order $2$, and since it interchanges the orbit $C$ and $Cb$ of $R(C)$, $x$ is semiregular.

Assume that $y$ is semiregular of order $8$. If $j$ is odd, then $i-ri+rj$ is odd and hence $y^2$ has order $8$ by Eq~(\ref{square2}), forcing $y$ has order $16$, a contradiction. Thus, $j$ is even, and hence $\g^{2(j+jr)}=1$. If $(i-ri+rj)/2$ is even, then $y^4=1$, contradicting $o(y)=8$. Thus, $(i-ri+rj)/2$ is odd. On the other hand, assume that  $j$ is even and $(i-ri+rj)/2$ is odd. By Eq~(\ref{square2}), $y^4=R(a)^4$, and hence $y$ has order $8$ and is semiregular.
\qed
\medskip

Recall that $H$ is regular on $\D_{16}$. Then we may let $H=\langle x,y\rangle$, where $o(x)=8$, $o(y)=2$ and $xx^y=1$. Since $H\leq A\leq \Hol(\D_{16})$, by Eq~(\ref{HolD_16C}) we may write $$x=x_1x_2 \mbox{ and } y=y_1y_2, \mbox{ where } x_1,y_1\in A1 \mbox{ and } x_2,y_2\in A2.$$ Since $A1\cong C_8^2$ and $A2\cong C_2^3$, we have $xx^y=(x_1x_2)(x_1x_2)^{y_1y_2}=x_1(x_1x_2)^{y_1y_2x_2}x_2=
x_1(x_1x_2^{y_1})^{y_2x_2}x_2=
x_1(x_1y_1^{-1}y_1^{x_2}x_2)^{y_2x_2}x_2=x_1(x_1y_1^{-1}y_1^{x_2})^{y_2x_2}=
x_1x_1^{x_2y_2}(y_1^{y_2x_2})^{-1}y_1^{y_2}$. Since $A1\unlhd \Hol(\D_{16})$,  $xx^y=1$ implies that
\begin{equation}\label{factorization}
x_1^{x_2}x_1^{y_2}=y_1^{-1}y_1^{x_2}
\end{equation}

Recall that $\{C,Cb\}$ is a complete imprimitive block system of $\Hol(\D_{16})$. Since $D$ fixes $1$, it fixes $C$ and $Cb$ setwise, but $R(b)$ interchanges them. By Eq~(\ref{HolD_16C}), $\Hol(\D_{16})=D\cup DR(b)$, and by the regularity of $H$, we only have two choices:

\medskip
\noindent {\bf Choice 1:} $x\in D$ and $y\in DR(b)$. By Claim~1,  $\bf x=R(a)^i\g^j\t$, where $i$ is odd and $j$ is even with $\t=1$ or $\b$, and  $\bf y=R(a)^k \g^\ell  \d R(b)$ for $\d\in B$, where $R(a)^\d=R(a)^r$ with $(1+r)\ell=0$ and $\ell=(1-r)k$.

\noindent {\bf Choice 2:} $x\in DR(b)$ and $y\in D$ (if $y\in DR(b)$ then replace $xy$ by $y$ as $H$ is dihedral). Note that $\{\lg a^2\rg,\lg a^2\rg a,\lg a^2\rg b,\lg a^2\rg ab\}$ is a complete imprimitive block system of $\Hol(\D_{16})$. By Claim~1, $\bf x=R(a)^i\g^j\d R(b)$, where $R(a)^\d=R(a)^r$, $j$ is even and $(i-ri+rj)/2$ is odd. By Eq~(\ref{square2}), $x^2$ has two orbits on $C$, that is, $\lg a^2\rg$ and $\lg a^2\rg a$, and $y$ interchanges the two orbits because $H$ is regular. Again by Claim~1, $\bf y=R(a)^k\g^\ell\a$, where $k$ is odd and $\ell$ is even.

\medskip
We finish the proof by considering four cases for $\Aut(\D_{16},S)$.

\medskip
\noindent {\bf Case 1:}  $\Aut(\D_{16},S)\leq \lg \g\rg$.

For Choice~1, $\g^j\t,\g^\ell  \d\in \Aut(\D_{16},S)$ as $R(\D_{16})\leq A$. It follows that $\t=1$ and $\d=1$, forcing $r=1$. Then $x=R(a)^i\g^j$, where $i$ is odd and $j$ is even, and $y=R(a)^k R(b)$ as $\ell=(1-r)k=0$. In particular, $x_1=R(a)^i\g^j$, $x_2=1$, $y_1=R(a)^k$ and $y_2=R(b)$.  Since $y_1^{-1}y_1^{x_2}=y_1^{-1}y_1=1$, Eq~(\ref{factorization}) implies that $1=x_1x_1^{R(b)}=(R(a)^i\g^j)(R(a)^i\g^j)^{R(b)}=R(a)^{j}\g^{2j}$.
It follows $j=0$ and $H=\lg R(a)^i, R(a)^k R(b)\rg=R(\D_{16})\unlhd A$, as required.

For Choice~2, $y=R(a)^k\g^\ell\a$, which is impossible because $\g^\ell\a\not\in \Aut(\D_{16},S)$.

\medskip
\noindent {\bf Case 2:}  $\Aut(\D_{16},S)\leq \langle \g^2, \g\a\b\rg$ or $\lg\g^2,\g\a\rg$.

Clearly, $\lg\g^2,\g\a\rg=\{1,\g^2,\g^4,\g^6,\g\a,\g^3\a,\g^5\a,\g^7\a\}$.
Recall that $\lg\g^2,\g\a\b\rg$ is the quaternion group of order $8$ and so
$\lg\g^2,\g\a\b\rg=\{1,\g^2,\g^4,\g^6,\g\a\b,\g^3\a\b,\g^5\a\b,\g^7\a\b\}$.

For Choice~1, $\g^j\t,\g^\ell  \d\in \Aut(\D_{16},S)$ implies that $\t=1$ and $\d=1$ as $\ell$ is even. It follows that $r=1$ and $\ell=(1-r)k=0$. Then a similar argument to Case~1 gives rise to  $H=R(\D_{16})\unlhd A$, as required. Also, Choice~2 cannot happen because $y=R(a)^k\g^\ell\a$ with $\ell$ is even and $\g^\ell\a\not\in \langle \g^2, \g\a\b\rg\cup \lg\g^2,\g\a\rg$.

\medskip
\noindent {\bf Case 3:}  $\Aut(\D_{16},S)$ is elementary abelian.

First we assume $\Aut(\D_{16},S)=\{1,\g^m\a,\g^{m+4}\a,\g^4\}$ for some odd $m$.

For Choice~1, $\g^j\t,\g^\ell  \d\in \Aut(\D_{16},S)$ implies that $\t=1$, $\g^j\in \lg\g^4\rg$ and $\d=1$ as $\ell$ is even, forcing $r=1$ and $\ell=0$. Then a similar argument to Case~1 gives rise to  $H=R(\D_{16})\unlhd A$, as required.
Also, Choice~2 cannot happen because  $y=R(a)^k\g^\ell\a$  with $\ell$ is even and $\g^\ell\a\not\in \Aut(\D_{16},S)$.

Now assume $\Aut(\D_{16},S)\not=\{1,\g^m\a,\g^{m+4}\a,\g^4\}$ for any odd $m$. Since $|\Aut(\D_{16},S)|\geq 4$, Observation implies that $i$ is even for every $\g^i\a^j\b^k\in \Aut(\D_{16},S)$.

Recall that $x\in H$ has order $8$. Since $x\in A=R(\D_{16})\rtimes \Aut(\D_{16},S)$ and $\Aut(\D_{16},S)$ is elementary abelian, $x^2$ is an element of order $4$ in $R(\D_{16})$, and hence $x^2=R(a^2)$ or $R(a^6)$. Write $L=\langle R(a^2)\rangle=\lg x^2\rg$. Note that $L\unlhd \Hol(\D_{16})$, $|R(\D_{16}):L|=4$ and $|L|=4$. We consider the quotient group $A/L$. For any $z\in \Hol(\D_{16})$  and any subgroup $T$ of  $\Hol(\D_{16})$, write $\overline{z}=zL$ and $\overline{T}=TL/L$. Then $\overline{H}=H/L\leq \overline{A}=A/L$. To prove $H\unlhd A$, it suffices to show that $\overline{A}$ is abelian.

Note that $\overline{R(\D_{16})}=\overline{R(C)}\times\overline{R(b)}\cong C_2^2$.
Since $R(g)^\d=R(g^\d)$ for any $R(g)\in R(\D_{16})$ and $\d\in \Aut(\D_{16})$, by definition, $\a$, $\b$ and $\g^2$ have trivial action on $\overline{R(\D_{16})}$, implying that $\Aut(\D_{16,S})$ has trivial action on $\overline{R(\D_{16})}$. It follows that $\overline{A}=\overline{R(\D_{16})}\times \overline{\Aut(\D_{16},S)}$, and since $\Aut(\D_{16},S)$ is elementary abelian, $\overline{A}$ is abelian, as required.

\medskip
\noindent {\bf Case 4:}  Others.

In this case, Cases~1-3 cannot happen. Then $\Aut(\D_{16},S)$ is not elementary abelian, and not a subgroup of $\lg\g\rg$, $\lg\g^2,\g\a\b\rg$ or $\lg\g^2,\g\a\rg$, that is, $\Aut(\D_{16},S)\nleqslant \lg\g\rg$, $\Aut(\D_{16},S)\nleqslant \lg\g^2,\g\a\b\rg$ and $\Aut(\D_{16},S)\nleqslant \lg\g^2,\g\a\rg$. We first prove the following claim.

\medskip
\noindent{\bf Claim 2:} One of the following holds: (1) $\{\g^4,\b\}\subseteq  \Aut(\D_{16},S)$; (2) $\{\g^2,\d\}\subseteq \Aut(\D_{16},S)$ with $1\not=\d\in  \langle \a,\b\rangle$; (3) $\g^i\b\in \Aut(\D_{16},S)$ for some odd $i$; (4) $\g^2\b\in\Aut(\D_{16},S)$.

\medskip
Let $\Aut(\D_{16},S)$ contain an element of order $8$. By Observation, $O_8=\{\g,\g^3,\g^5$, $\g^7$, $\g\b$, $\g^3\b$, $\g^5\b$, $\g^7\b\}$. Since $\Aut(\D_{16},S)\nleqslant \lg\g\rg$, we have $\g^i\b\in \Aut(\D_{16},S)$ for some odd $i$, and Claim~2~(3) follows. Note that $(\g^2\b)^{-1}=\g^6\b$. If $\g^2\b$ or $\g^6\b\in \Aut(\D_{16},S)$, then $\g^2\b\in \Aut(\D_{16},S)$ and Claim~2~(4) follows. By Observation, $O_4=\{\g^2, \g^6, \g^2\b, \g^6\b, \g\a\b, \g^3\a\b, \g^5\a\b, \g^7\a\b\}$. Since $\Aut(\D_{16},S)$ is not elementary abelian, we  assume that $\Aut(\D_{16},S)$ has no element of order $8$, $\g^2\b,\g^6\b\not\in \Aut(\D_{16},S)$, and either $\g^2\in \Aut(\D_{16},S)$ or $\g^j\a\b\in \Aut(\D_{16},S)$ for some odd $j$.

Let $\g^j\a\b\in \Aut(\D_{16},S)$ for some odd $j$.  Since $\Aut(\D_{16},S)\nleqslant \lg\g^2,\g\a\b\rg=\{1,\g^2,\g^4$, $\g^6$, $\g\a\b$, $\g^3\a\b$, $\g^5\a\b$, $\g^7\a\b\}$, $\Aut(\D_{16},S)$ has an involution $\theta\not=\g^4$.  By Observation, $O_2=\lg\g\rg\a\cup \lg\g^2\rg\a\b\cup \{\b,\g^4,\g^4\b\}$.
If $\theta\in\{\b,\g^4\b\}$, then $\{\g^4,\b\}\subseteq \Aut(\D_{16},S)$ and Claim~2~(1) follows. If $\theta\in \lg\g^2\rg\a\b$, then $\g^j\a\b\theta\in \{\g,\g^3,\g^5,\g^7\}$ and $\g^i\in \Aut(\D_{16},S)$ for some odd $i$, which is impossible as $o(\g^i)=8$. If  $\theta\in \lg\g\rg\a$, then $\theta\g^j\a\b\in \lg\g\rg\b\subseteq \Aut(\D_{16},S)$. Since $\g^2\b,\g^6\b,\g^i\b\not\in \Aut(\D_{16},S)$ for any odd $i$ ($o(\g^i\b)=8$), we have $\g^4\b\in \Aut(\D_{16},S)$, and $\{\g^4,\b\}\subseteq \Aut(\D_{16},S)$. The Claim~2~(1) follows.

Assume $\g^i\a\b\not\in \Aut(\D_{16},S)$ for all odd $i$. Then $\g^2\in\Aut(\D_{16},S)$ and $\Aut(\D_{16},S)$ has only two elements of order $4$, that is, $\g^2$ and $\g^6$. Since $\Aut(\D_{16},S)\nleqslant \lg\g\rg$, $\Aut(\D_{16},S)$ has an involution $\rho\not=\g^4$. Recall that $O_2=\lg\g\rg\a\cup \lg\g^2\rg\a\b\cup \{\b,\g^4,\g^4\b\}$. If $\rho\in \lg\g^2\rg\a\cup\lg\g^2\rg\a\b\cup \{\b,\g^4\b\}$ then $\{\g^2,\d\}\subseteq \Aut(\D_{16},S)$ with $\d=\a,\a\b,\b$, and  Claim~2~(2) follows. If $\rho\in \{\g\a,\g^3\a,\g^5\a,\g^7\a\}$, then also some involution in $\lg\g^2\rg\a\cup\lg\g^2\rg\a\b\cup \{\b,\g^4\b\}$ belongs to $\Aut(\D_{16},S)$  as $\Aut(\D_{16},S)\nleqslant\lg\g^2,\g\a\rg$, and again Claim~2~(2) follows.  \qed

\medskip
For Claim~2~(1),  $\{\g^4,\b\}\subseteq  \Aut(\D_{16},S)$. Write $K=\langle a^4\rg$. Then $F_{\D_{16}}(\g^4)=C$, where $F_{\D_{16}}(\g^4)$ is the fixed point set of $\g^4$ in $\D_{16}$, as defined in Proposition~\ref{MainPro}. Clearly, for every $g\in \D_{16}$, $\lg \g^4\rg$ fixes $Kg$ pointwise or has $Kg$ as an orbit, and $\b$ fixes $Kg$ setwise. By Proposition~\ref{MainPro}~(3), $\Ga$ is non-normal, a contradiction. For Claim~2~(2), $\{\g^2,\d\}\subseteq \Aut(\D_{16},S)$ with $1\not=\d\in  \langle \a,\b\rangle$, and by setting $K=\langle a^2\rg$, we have the same contradiction as above.

Let $T_1\subseteq C$ and $T_2\subseteq Cb$. Denote by $[T_1,T_2]$ the induced subdigraph of $\Ga$, which has vertex set $T_1\cup T_2$ and all arcs of $\Ga$ between $T_1$ and $T_2$.  Clearly, these arcs corresponds to involutions in $Cb$, and hence $[T_1,T_2]$ is a graph.

For Claim~2~(3), $\g^i\b\in \Aut(\D_{16},S)$ with $i$ odd. Then $\g^i\b$ has order $8$. Since $\g$ fixes $C$ pointwise, $\g^i\b$ and $\b$ have the same action on $C$, and $\lg\g^i\b\rg$ has orbits of length $1$ or $2$ on $C$. Since $o(\g^i\b)=8$, $\lg\g^i\b\rg$ is transitive on $Cb$, implying that
either $1$ is adjacent to every vertex in $Cb$ or has no neighbours in $Cb$. Since $R(C)\leq A$,  $[C,Cb]$ is an empty graph or $K_{8,8}$, the complete bipartite graph of order $16$.  Thus, the restriction of $\g^i\b$ on $C$ can be extended to an automorphism of $\Ga$ by fixing $Cb$ pointwise, say $\overline{\b}$. Then  for every $c\in C$, we have
$c^{\overline{\b}}=c^{\g^i\b}=c^{\b}$ and $(cb)^{\overline{\b}}=cb$. Since $\b\not=1$, $\overline{\b}$ has order $2$. On the other hand, since $\overline{\b}$ fixes $1$ and $\Ga$ is normal, we have $\overline{\b}\leq \Aut(\D_{16},S)$. Since $\overline{\b}$ fixes $Cb$ pointwise, it fixes $\D_{16}=\lg Cb\rg$ pointwise, and hence $\overline{\b}=1$, a contradiction.

For Claim~2~(4), $\g^2\b\in\Aut(\D_{16},S)$. Then $\g^2\b$ has order $4$. Write $L=\lg a^2\rg$. Note that $\lg \g^2\b\rg$ fixes the complete imprimitive block system $\{L,La,Lb,Lab\}$ pointwise, and is transitive on both $Lb$ and $Lab$. Then each of the induced subdigraphs $[L,Lb]$, $[L,Lab]$, $[La,Lb]$ and $[La,Lab]$ is $K_{4,4}$ or an empty graph. Then the restriction of $\g^2\b$ on $C$ can be extended to an automorphism of order $2$ of $\Ga$ by fixing $Cb$ pointwise, also denoted by $\overline{\b}$, that is,
$c^{\overline{\b}}=c^{\g^2\b}=c^{\b}$ and $(cb)^{\overline{\b}}=cb$ for every $c\in C$. Then a similar argument to the above paragraph gives rise to the contradiction $\overline{\b}=1$. This completes the proof. \qed

\bigskip
\noindent{\bf Proof of Theorem~\ref{mainth2}:} Let $\D_{2n}$ be a dihedral group of order $2n$. Let $n\ge 6$ be even and $n\neq 8$. By Lemma~\ref{mainlem}, $\D_{2n}$ is an NNN-group, and hence an NNND-group because NNN-group implies NNND-group. On the other hand, assume that  $\D_{2n}$ is an NNND-group. By Lemma~\ref{lem2}, $n\not=8$. Furthermore, $\D_{2n}$ is not an NDCI-group, and by Proposition~\ref{NDCID}, $n$ is even but $n\not=2,4$. It follows that $n\ge 6$ is even and $n\neq 8$. Clearly, this is also true if $\D_{2n}$ is an NNN-group. \qed

\section*{Acknowledgments}
This work was partially supported by the National Natural Science Foundation of China (12331013, 12271024, 12425111, 12301461) and the 111 Project of China (B16002).

\end{document}